\documentclass[journal,twoside,web]{ieeecolor}
\usepackage{generic}
\usepackage{cite}
\usepackage{amsmath,amssymb,amsfonts}
\usepackage{graphicx}
\usepackage{textcomp}
\usepackage{bbm}
\usepackage[dvipsnames]{xcolor}
\usepackage{setspace} 
\usepackage{graphicx}        
\usepackage{amssymb}
\usepackage{graphicx}
\usepackage{algpseudocode}
\usepackage{algorithm}
\usepackage{dsfont}
\usepackage{dblfloatfix} 

\usepackage{amsmath,todonotes}
\usepackage{mathtools}
\usepackage{comment}
\usepackage{multirow}
\usepackage{booktabs}
\usepackage{stmaryrd}
\usepackage[hidelinks]{hyperref}
\usepackage{pgfplots}

\usetikzlibrary{calc}
\usetikzlibrary{shapes.misc}
\tikzset{cross/.style={cross out, draw=black, minimum size=2*(#1-\pgflinewidth), inner sep=0pt, outer sep=0pt},cross/.default={1pt}}

\usepackage{array}

\newcolumntype{L}[1]{>{\raggedright\let\newline\\\arraybackslash\hspace{0pt}}m{#1}}
\newcolumntype{C}[1]{>{\centering\let\newline\\\arraybackslash\hspace{0pt}}m{#1}}
\newcolumntype{R}[1]{>{\raggedleft\let\newline\\\arraybackslash\hspace{0pt}}m{#1}}

\let\proof\relax
\let\endproof\relax

\usepackage{amsthm}

\theoremstyle{plain}

\newtheorem{theorem}{Theorem}
\newtheorem{lemma}{Lemma}
\newtheorem{proposition}{Proposition}

\newtheorem{examplex}{Example}
\newenvironment{example}
  {\pushQED{\qed}\examplex}
  {\popQED\endexamplex}

\newcommand\res{r}
\newcommand\A{\mathcal{A}}
\newcommand\U{\mathcal{U}}
\newcommand\X{\mathcal{X}}

\newcommand\T{\mathcal{T}}
\newcommand\Q{\overline{\mathcal{M}}}

\newcommand\M{\mathcal{M}}
\newcommand\supp{\operatorname{supp}}

\let\ForAll\forall
\renewcommand\forall{\ForAll\,}

\def\BibTeX{{\rm B\kern-.05em{\sc i\kern-.025em b}\kern-.08em
    T\kern-.1667em\lower.7ex\hbox{E}\kern-.125emX}}
\begin{document}

\title{Strategic Monitoring of Networked Systems with Heterogeneous Security Levels}

\author{\makebox[\linewidth][c]{Jezdimir Milo\v{s}evi\'{c}, Mathieu Dahan, Saurabh Amin, and Henrik Sandberg}
\thanks{Manuscript received April 6, 2023. The work of JM and HS was supported in part by the Swedish Civil Contingencies Agency (project CERCES2) and Digital Futures (project DEMOCRITUS). The work of MD was supported by the Georgia Tech new faculty start up grant. The work of SA was supported in part by AFOSR grant ``Building Attack Resilience into Complex Networks'' and NSF grant \# 2039771.}
\thanks{Jezdimir Milo\v{s}evi\'{c} was with the School of Electrical Engineering and Computer Science, KTH Royal Institute of Technology, Stockholm, Sweden. He is now with Scania Autonomous Transport Solutions (e-mail:  jezdimirmilosevic@gmail.com).}
\thanks{Mathieu Dahan is with the School of Industrial and Systems Engineering, Georgia Institute of Technology, Atlanta, GA 30332 USA (e-mail: mathieu.dahan@isye.gatech.edu).}
\thanks{Saurabh Amin is with the Department of Civil and Environmental Engineering, Massachusetts Institute of Technology, Cambridge, MA 02139 USA (e-mail: amins@mit.edu).}
\thanks{Henrik Sandberg is with the School of Electrical Engineering and Computer Science, KTH Royal Institute of Technology, Stockholm, Sweden (e-mail: hsan@kth.se).}}

\maketitle

\begin{abstract} 
\textbf{We consider a strategic network monitoring problem involving the operator of a networked system and an attacker. The operator aims to randomize the placement of multiple protected sensors to monitor and protect components that are vulnerable to attacks. We account for the heterogeneity in the components' security levels and formulate a large-scale maximin optimization problem. After analyzing its structure, we propose a three-step approach to approximately solve the problem. First, we solve a generalized covering set problem and run a combinatorial algorithm to compute an approximate solution. Then, we compute approximation bounds by solving a nonlinear set packing problem.
To evaluate our solution approach, we implement two classical solution methods based on column generation and multiplicative weights updates, and test them on real-world water distribution and power systems. Our numerical analysis shows that our solution method outperforms the classical methods on large-scale networks, as it efficiently generates solutions that achieve a close to optimal performance and that are simple to implement in practice.}
%
%
%
%
\end{abstract}

 
\section{Introduction}

\IEEEPARstart{L}{arge-scale} networked systems, such as water distribution or power systems, are lucrative targets for~malicious~attackers. 
Indeed, several attacks have already occurred~\cite{slay2007lessons,CSIS2022}, demonstrating the need for rapid development of effective defense strategies  for these systems. 
 An essential component of every defense strategy is attack detection~\cite{nist}, which can be accomplished by allocating multiple sensors for system monitoring. 
Yet, for large-scale systems, the number of available sensors is likely to be insufficient to ensure complete coverage.
A strategic adversary can exploit this deficiency to target the most critical components of the system that are unmonitored.  While the literature on sensor placement has studied fixed \cite{stack_metju,pirani2018game} and randomized sensing strategies \cite{2017arXiv170500349D,krause2011randomized}, none of these works has accounted for the heterogeneous criticality of the network components, and its impact on monitoring decisions. Hence, the question that arises is how to strategically allocate multiple sensors to monitor a large-scale system of critical networked components.
%

Specifically, we focus on determining a randomized scheduling of monitoring operations to maximize the protection of a networked system against a strategic adversary. In this network monitoring problem, the network operator seeks to randomize the placement of multiple protected sensors to monitor the components and protect them against an attacker who aims to sabotage a system component while remaining undetected. To better account for the attacker's behavior and preference, we consider that the components have possibly heterogeneous security levels, which represent the corresponding levels of effort required by an adversary to target them. This model is motivated by risk management, where one first identifies critical system components in a risk assessment, and then allocates protection resources accordingly \cite{nist}. By monitoring a component with a protected sensor, the operator maximizes the security level of that component. Then, to protect the system against a strategic adversary, we assume that the operator aims to maximize the lowest expected security level across all components.

 \subsubsection*{Our contributions}   
In our first set of contributions, we analyze the structure of the network monitoring problem and gather insights on the impact of the problem characteristics on the operator's optimal strategy. This leads us to a three-step approach for efficiently computing an approximate monitoring strategy and evaluating its optimality gap. First, we formulate and solve a mixed-integer program (MIP) to compute generalized covering sets and obtain marginal probabilities of placing sensors at each location. Second, we employ the combinatorial algorithm of \cite{dziubinski2018hide} to construct a randomized monitoring strategy that is consistent with the marginal probabilities computed from the MIP. Finally, by leveraging our structural insights, we formulate a MIP to compute nonlinear set packings and obtain an optimality gap on our approximate solution to the network monitoring problem.

In our second set of contributions, we refine our solution approach and analysis for problem instances with special structures. In particular, we find that our solution approach optimally solves the network monitoring problem when the sets of components that are monitored from each sensor location are mutually disjoint. Furthermore, when security levels are identical, our solution approach can be simplified by solving the minimum set cover and maximum set packing problems.

We also describe and implement two classical solution methods for solving the network monitoring problem based on column generation and multiplicative weights updates. In particular, we show that these algorithms can be utilized, as the pricing problem for the column generation algorithm and the best-response problem for the multiplicative weights update algorithm, can both be formulated as maximum weighted covering problems.


In our last set of contributions, we evaluate our solution approach based on generalized covering sets and nonlinear set packings and compare it with the two classical solution methods using two security applications. In the first application, we consider the problem of detecting contaminants injected in a water distribution system. In the second application, we consider the problem of protecting actuators in a large-scale power system against extended replay attacks, during which the attacker uses actuators in addition to sensors to lead an undetectable attack~\cite{mo2015physical,milovsevic2020security}. Our computational results show that our solution approach provides excellent solutions to the network monitoring problem and outperforms the two classical methods in different aspects: It generates solutions with low optimality gaps, and its running time is marginally impacted by the number of available sensors. In contrast, the running time of the other methods exponentially increases with respect to the number of sensors. Finally, the monitoring strategies generated by our solution approach are significantly simpler to implement in practice compared to the ones generated by the classical methods, thus reinforcing the value of our solution approach derived from structural insights.



%
%

 \subsubsection*{Related work}   
 
 
Optimization and game theoretic models have been developed for studying various security related problems,
including the development of defense strategies~\cite{zhu2015game,BaykalGursoy2014469,Goyal:2014aa,doi:10.1287/opre.1070.0434}, the design of anomaly detectors~\cite{miao2018hybrid}, the allocation of security budget~\cite{hota2016optimal},  and the placement of sensors~\cite{stack_metju,pirani2018game,2017arXiv170500349D,krause2011randomized}. 
For instance, \cite{doi:10.1287/inte.1060.0252,Bier_2011,doi:10.1111/risa.12333,doi:10.1287/trsc.2017.0749} consider bilevel and trilevel optimization problems to model defender-attacker interactions where each player selects a pure strategy. In contrast, we focus on randomized strategies, which are recognized to be more effective when the number of sensors to deploy is limited~\cite{zhu2015game, hota2016optimal}.

Other models that have been studied include inspection problems \cite{doi:10.1287/opre.43.2.243,doi:10.1287/opre.46.2.184,smith2008algorithms} and search games \cite{Gal20140062,doi:10.1287/opre.2019.1853}, in which a searcher is concerned with the optimal way of looking for a hidden adversary in a search space; see for example \cite{lidbetter2013search,clarkson2022classical,alpern2006theory,hohzaki2016search}.  In particular, the recent literature on hide-and-seek games investigates the problem of coordinating multiple inspection resources for detecting hidden objects in boxes \cite{dziubinski2018hide,gal2014succession, https://doi.org/10.48550/arxiv.2301.11541}. However, such problems do not account for the impact of sensing range and the network topology on the monitoring strategies. Still, we leverage the algorithm of \cite{dziubinski2018hide} as part of our solution approach for computing randomized monitoring strategies with small supports.

Finally, the existing literature on sensor placement considers developing both fixed \cite{stack_metju,pirani2018game} and randomized monitoring strategies \cite{2017arXiv170500349D,krause2011randomized}.
The placement problems considered previously aimed to minimize estimation error~\cite{7524914}, achieve optimal coverage~\cite{cortes2004coverage}, detect faults~\cite{perelman2016sensor}, or improve the system's security~\cite{TAC_2019_Jezdimir}.  
The study~\cite{TAC_2019_Jezdimir} is related to ours, and proposes two static sensor placement strategies to improve the security of actuators in large-scale systems. 
Our focus is instead on randomized strategies, which can be translated into random scheduling of inspections that can be performed on a day-to-day basis by utility or security personnel \cite{pita2008deployed}. 

Our problem is also related to the simultaneous zero-sum game in~\cite{2017arXiv170500349D}, where the operator places multiple sensors to maximize the number of detected (homogeneous) attacks. Recent work extended this model by considering the imperfect detection by sensors  \cite{9867753}, and the different types of sensing resources \cite{McCann2022}, for which they derive heuristic approaches. In addition, \cite{pirani2021strategic} formulated a game in which sensors are positioned in the nodes of a networked control system to detect attacks on them, and considered imperfect detection through a linear filter that processes the sensors' measurements to detect attacks. The authors derived equilibrium results using tools from structured systems and graph theory. However, to the best of our knowledge, no article investigates the coordination of \emph{multiple} sensing resources to \emph{strategically} monitor a \emph{networked} system with components of \emph{heterogeneous} criticalities. 

Preliminary results (without proofs) have been presented in \cite{jezd_cdc}. Specifically, \cite{jezd_cdc} investigate a simpler solution approach based on minimal set covers, maximal set packings, and column generation. They also provide approximation bounds under restrictive assumptions. By further leveraging the problem structure, we generalize their results and obtain a more efficient solution approach with stronger approximation bounds. We also provide a more comprehensive evaluation of our solution approach through an extended computational study that involves two other solution algorithms and two security applications.

 \subsubsection*{Organization} 
Section~\ref{section:game} introduces and formulates the strategic network monitoring problem. In Section \ref{sec:Method}, we derive and analyze our solution approach based on generalized covering sets and nonlinear set packings. We then present two alternative solution methods based on column generation and multiplicative weights updates in Section \ref{sec:Alt}. Section \ref{section:computations} evaluates and compares these solution approaches through a computational study. Finally, Section~\ref{section:ch8_conclusion} summarizes our findings and provides avenues for future work.


 \section{Problem Formulation}\label{section:game}
 
 We consider a networked system consisting of a set of components $\mathcal{U}$. Each component $u \in \U$ faces a risk of attack by an adversary who wishes to disrupt or take control of the network. For every $u \in \U$, let $\varphi_u \in [0,1]$ denote its \emph{security level}, i.e., the effort that the attacker needs to spend in order to successfully target $u$ and compromise its operation. Security levels can for example be estimated based on previously deployed security measures, configuration of the system, or monetary cost for conducting an attack against the components. A component $u$ for which $\varphi_u = 0$ (resp. $\varphi_u = 1$) is unsafe (resp. highly secure). 

The network operator can increase the security levels of some components by positioning $\res$ sensors within a set of locations $\X$. Specifically, a sensor positioned at location $x \in \X$ can monitor and secure a subset of components $U_x \subseteq  \mathcal{U}$. 
For every $x \in \X$, we refer to $U_x$ as the monitoring set of  $x$. Typically, monitoring sets are computed from the characteristics of the sensors and the network. If multiple sensors are positioned at locations $X \subseteq \mathcal{X}$, then the monitored and secured components are given by $U_X \coloneqq \cup_{x \in X}U_x$. 

We let $\A \coloneqq \{X \subseteq \X \ | \ |X| \leq \res\}$ denote the set of actions of the operator. Thus, for every sensor positioning $X \in \A$ and every component $u \in \U$, the resulting \emph{post-security level} of $u$ is given by $f(X,u) \coloneqq \varphi_u \mathds{1}_{\{u \notin U_X\}} + \mathds{1}_{\{u \in U_X\}}$. Specifically, the post-security levels of all monitored components are maximal and equal to one. 
Without loss of generality, we assume the following: 1) $U_x \neq  \emptyset$ for every $x \in \X$; 2) every component can be monitored from at least one sensor location; and 3) $\varphi_u < 1$ for every $u \in \U$.

Next, we illustrate the monitoring model with an example.
\begin{example}
Consider the networked system and monitoring model illustrated in Fig. \ref{fig:example}.

\begin{figure}[htbp]
    \centering
        \begin{tikzpicture}[scale=1.25,main node/.style={circle,draw,inner sep = 0.15cm},main node2/.style={circle,draw,inner sep = 0.10cm},]
            \fill[fill=gray!20!white] (-0.5, 0) ellipse [x radius = 15mm, y radius = 6mm];
            \fill[fill=gray!20!white] (1.5, 0) ellipse [x radius = 15mm, y radius = 6mm];
            \fill[fill=gray!20!white] (1.5, -1.3) ellipse [x radius = 6mm, y radius = 15mm];
            
            \draw (-0.5, 0) ellipse [x radius = 15mm, y radius = 6mm];
            \draw (1.5, 0) ellipse [x radius = 15mm, y radius = 6mm];
            \draw (1.5, -1.3) ellipse [x radius = 6mm, y radius = 15mm];
            
            
            \node[scale =0.9] (1) at (-0.5, 0.4) {$x_1$};
            \node[scale =0.9] (2) at (2.5, 0.15) {$x_2$};
            \node[scale =0.9] (3) at (1.8, -1) {$x_3$};
            
            \node[circle, draw, inner sep = 0.07cm, fill=blue!50!white, scale=0.9] (e_1) at (-1.5, 0) {$u_1$};
            \node[circle, draw, inner sep = 0.07cm, fill=Green!90!white, scale=0.9] (e_2) at (-0.5, -0.3) {$u_2$};
            \node[circle, draw, inner sep = 0.07cm, fill=red!90!white, scale=0.9] (e_3) at (0.5, 0) {$u_3$};
            \node[circle, draw, inner sep = 0.07cm, fill=red!90!white, scale=0.9] (e_4) at (1.5, -0.3) {$u_4$};
            \node[circle, draw, inner sep = 0.07cm, fill=blue!50!white, scale=0.9] (e_5) at (2.0, 0.25) {$u_5$};
            \node[circle, draw, inner sep = 0.07cm, fill=Green!90!white, scale=0.9] (e_6) at (1.2, -1.3) {$u_6$};
            \node[circle, draw, inner sep = 0.07cm, fill=blue!50!white, scale=0.9] (e_7) at (1.5, -2.3) {$u_{7}$};
            
            \node (green sensor) at (1.8, -1.3) {\includegraphics[scale=0.08]{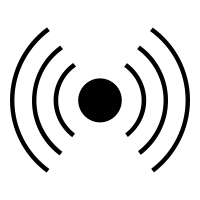}};
            \node (yellow sensor) at (2.5, -0.15) {\includegraphics[scale=0.08]{Sensor.png}};            

        \end{tikzpicture}
    \caption{Networked system containing 3 sensor locations and 7 components.} 
    \label{fig:example}
\end{figure}

In this example, the set of sensor locations is $\X = \{x_1,x_2,x_3\}$, the set of components is $\U = \{u_1,\dots,u_7\}$, and the monitoring sets are $U_{x_1} = \{u_1,u_2,u_3\}$, $U_{x_2} = \{u_3,u_4,u_5\}$, and $U_{x_3} = \{u_4,u_6,u_7\}$. The security levels are given by $\varphi_{u_3} = \varphi_{u_4} = 0.2$ (red in the figure), $\varphi_{u_1} = \varphi_{u_5} = \varphi_{u_7} = 0.5$ (blue in the figure), and $\varphi_{u_2} = \varphi_{u_6} = 0.8$ (green in the figure).
We assume that the operator has two sensors that they place on $X = \{x_2,x_3\}$. Then, each component $u$ in $U_X = \{u_3,u_4,u_5,u_6,u_7\}$ is monitored and has a post-security level equal to $f(X,u)= 1$. On the other hand, the remaining components $u_1$ and $u_2$ are not monitored, and their post-security levels are equal to their original security levels: $f(X,u_1) = 0.5$ and $f(X,u_2) = 0.8$. Note that $u_1$ has the lowest resulting post-security level and is the most vulnerable to attacks. 
%
%
%
%
\end{example}

In such combinatorial security settings, the network operator can significantly benefit from randomizing the positioning of their sensors \cite{washburn1995two, pita2008deployed,  gupta2016dynamic}. Thus, we consider the problem where the operator is interested in positioning their sensors in a randomized manner to protect the network against an attacker who would target the most vulnerable component. In other words, the operator aims to maximize the lowest expected post-security level across all network's components. Let $\Delta^1\coloneqq\{\sigma^{1} \in [0,1]^{|\mathcal{A}|}  
\ | \
 \sum_{X \in \mathcal{A}} \sigma^1_{X} =1 \}$ denote the set of mixed strategies of the operator, i.e., the set of probability distributions over $\A$; here, $\sigma^1_X$ denotes the probability that the operator positions sensors on locations $X \in \A$. Then, the problem of strategic network monitoring can be formulated as
 \begin{align}
 \label{eq:P}
\max_{\sigma^1 \in \Delta^1} \min_{u \in \U} \mathbb{E}_{X\sim\sigma^1}[f(X,u)]. \tag{$\mathcal{M}$}
 \end{align}


Problem \eqref{eq:P} can be formulated as a linear program with $|\A| +1$ variables and $|\U|+1$ constraints. However, the cardinality of $\A$ grows combinatorially with respect to the number of  available sensors. 
Thus, standard methods for solving linear programs cannot be used for large-scale systems, and different approaches are needed to compute an optimal or approximate monitoring strategy.
 In Section~\ref{sec:Method}, we develop an approach based on mixed-integer programs (MIPs) to compute an approximate solution of \eqref{eq:P}, while in Section~\ref{sec:Alt} we discuss alternative approaches based on column generation and multiplicative weights updates.

%

\section{Main Results}  \label{sec:Method}

In this section, we analyze the structure of the network monitoring problem \eqref{eq:P} and derive an approximate solution approach. Our approach first solves a generalized covering set problem and runs a combinatorial algorithm to construct a monitoring strategy with small support. Then, it computes an optimality gap by solving a nonlinear set packing problem. We further investigate problem instances with special structures, for which we refine our approach and optimality gaps.

%

\subsection{Preliminaries}

To carry out our analysis, we first introduce additional quantities and notations:
For any optimization problem $(\mathcal{O})$, we denote by $z^*_{(\mathcal{O})}$ its optimal value. Then, for every location $x \in \X$, let $\overline{\varphi}_x \coloneqq \min_{u \in U_x} \varphi_u$ denote the lowest security level of a component in the monitoring set of $x$. This quantity will impact the operator's monitoring strategy, as the objective is to maximize the worst-case component post-security level. We quantify the \emph{criticality} of a location $x \in \X$ (resp. component $u \in \U$) from the value of $\overline{\varphi}_x$ (resp. $\varphi_u$): If $\overline{\varphi}_x$ (resp. $\varphi_u$) is low, we say that  $x \in \X$ (resp. $u \in \U$) is more critical. We additionally define, for every sensor positioning $X \in \A$ and every location $x \in \X$, $\bar{f}(X,x) \coloneqq \overline{\varphi}_x \mathds{1}_{\{x \notin X\}} + \mathds{1}_{\{x \in X\}}$ to be the lowest post-security level within the monitoring set $U_x$.

For any subset of locations $C \subseteq \X$ (resp. any subset of components $T \subseteq \U$), we denote as $\overline{\varphi}_C$ (resp. $\varphi_T $) the following vector $(\overline{\varphi}_x)_{x \in C}$ (resp. $(\varphi_u)_{u \in T}$), sorted in nondecreasing order.
For every subset of components $T\subseteq \U$, we also define the following quantity $S_T \coloneqq \sum_{u \in T}(1-\varphi_u)^{-1}$.


For any vector $y \in \mathbb{R}^n$, we denote its support as $\supp(y) \coloneqq \{i \in \llbracket1,n\rrbracket \ | \ y_i \neq 0\}$, where $\llbracket1,n\rrbracket = \{1,\dots,n\}$.
Then, given a monitoring strategy $\sigma^1 \in \Delta^1$, we define its \emph{node basis} as $\X(\sigma^1) \coloneqq \cup_{X \in \supp(\sigma^1)} X$, that is, the set of nodes that are monitored with positive probability under the monitoring strategy $\sigma^1$. Note that node bases provide a good indicator of the practical implementability of a monitoring strategy: A strategy that monitors a small number of locations can easily be implemented in practice via a randomized scheduling of operations \cite{pita2008deployed}.

We consider the collection of set covers $\mathcal{S}$ of our networked system. Specifically, a subset of locations $C \subseteq \X$ is a \textit{set cover} if $U_C = \U$, i.e., every component is monitored from at least one location in $C$. Let $n^*$ denote the minimum set cover size. We note that if $r \geq n^*$, then an optimal solution of the network monitoring problem \eqref{eq:P} consists of placing sensors on a set cover and ensuring a post-security level of 1 for every component. Henceforth, we assume that $r < n^*$.

Finally, we consider the collection of set packings $\T$ of our networked system. Specifically, a subset of components $T \subseteq \U$ is a \textit{set packing} if $|T \cap U_x|\leq 1$ holds for every $x \in \mathcal{X}$. Set packings indicate the \emph{spread} of the network, which directly impacts the monitoring strategy of the operator: Since components in a set packing must be monitored from distinct sensor locations, a network with large set packings is more challenging to protect. We denote the maximum set packing size as $m^*$.
%
%
%
%

\subsection{Structural Analysis}

We first analyze the structure of the network monitoring problem \eqref{eq:P} and derive bounds on its optimal value.
To obtain a lower bound on $z^*_{(\M)}$, we simplify the problem by supposing that when the operator places a sensor at location $x \in \X$, they protect the component $u$ within the monitoring set $U_x$ that is the most critical. The premise is that the post-security level achieved for any component in $U_x$ will be at least that of component $u$. Formally, we relax the objective function of \eqref{eq:P} and consider the following problem:
%
%
%
%
%
%
\begin{align}
\label{eq:Q}
 \max_{\sigma^1\in \Delta^1} \min\bigg\{\ \min_{\quad \mathclap{x \in \X(\sigma^1)}\quad } \mathbb{E}_{X\sim\sigma^1}[\bar{f}(X,x)], \ \min_{\quad \mathclap{u \notin U_{\X(\sigma^1)}}\quad }\varphi_u\bigg\}.
\tag{$\overline{\mathcal{M}}$}
\end{align}

Problem \eqref{eq:Q} selects a monitoring strategy $\sigma^1$ that partitions the set of components $\U$ into two subsets: the subset of components $\U\setminus U_{\X(\sigma^1)}$ that are never monitored by $\sigma^1$ and for which the post-security levels are given by the corresponding original security levels; and the subset of components $U_{\X(\sigma^1)}$ that are monitored with positive probability. In this latter set, $\sigma^1$ guarantees a post-security level for every component in $U_{\X(\sigma^1)}$ by improving the lowest security level within each monitoring set $U_x$ for $x \in \X(\sigma^1)$.

In fact, we provide in the next lemma a different perspective on problem \eqref{eq:Q} to guide our analysis:
\begin{lemma}
\label{lem:Angle}
Problem \eqref{eq:Q} satisfies the following property:
\begin{align*}
z^*_{(\Q)}  = \max_{C \subseteq \X} \min\Big\{z^*_{(\mathcal{R}_C)},\min_{u \notin U_C} \varphi_u\Big\},
\end{align*}

where for every $C \subseteq \X$, $z^*_{(\mathcal{R}_C)}$ is the optimal value of
\begin{align}
\label{eq:R}
\max_{\sigma^1 \in \Delta^1} \min_{x \in C} \mathbb{E}_{X \sim \sigma^1}[\bar{f}(X,x)].\tag{$\mathcal{R}_C$}
\end{align}

\end{lemma}

%



Thus, \eqref{eq:Q} can be solved by first determining the set of locations $C$ over which to randomize the placement of sensors, and then determining the monitoring strategy that optimizes \eqref{eq:R}. Given a preselected subset of $n$ locations $C  =\{x_1,\dots,x_n\} \subseteq \X$, indexed such that $\overline{\varphi}_{x_1} \leq \cdots \leq \overline{\varphi}_{x_n}$,
%
%
%
%
%
we provide some intuition on the monitoring strategy $\sigma^1$ that maximizes $\min_{x \in C} \mathbb{E}_{X \sim \sigma^1}[\bar{f}(X,x)]$. Initially, the goal is to randomize the placement of the sensors among the locations in $C$ so as to equalize the post-security level of the most critical component in each of their monitoring sets.
However, because of the heterogeneity of the security levels, the operator should instead focus on the most critical locations and not monitor the most secure locations. The number of most critical locations in $C$  to monitor is given by 
\begin{align*}
k^*_C\coloneqq \max\bigg\{k \in \llbracket 1,n\rrbracket \ \bigg| \ \overline{\varphi}_{x_k} \leq 1 - \frac{k - \res}{\sum_{l=1}^k (1-\overline{\varphi}_{x_l})^{-1}}\bigg\}.
\end{align*}
For simplicity, we also denote $S^{k^*}_C \coloneqq \sum_{l=1}^{k^*_C}(1-\overline{\varphi}_C)^{-1}$.

%



Using this property and the collection of set packings $\mathcal{T}$, we can then derive the following main theorem: 
 \begin{theorem}
 \label{prop:bounds}
The optimal value of the network monitoring problem \eqref{eq:P} is bounded as follows:
 \begin{align*}
 z_{(\Q)}^* &= \max_{C \subseteq \X} \min\bigg\{1 - \frac{k^*_C-\res}{S^{k^*}_C},\min_{\mathclap{u \notin U_C}}\varphi_u\bigg\} \\
& \leq z_{(\M)}^*\leq \min\bigg\{1,\min_{T \in \T} 1 -  \frac{|T|-\res}{S_T}\bigg\}.
 \end{align*}
 
 Furthermore, for any set of $n > r$ locations $C = \{x_1,\dots,x_n\}\subseteq \X$ indexed in nondecreasing order of their criticalities, any monitoring strategy $\sigma^{1^\prime}$ satisfying the following conditions is an optimal solution of \eqref{eq:R}:
 \begin{align}
\mathbb{P}_{\hspace{-0.05cm}X\sim\sigma^{1^\prime}}(x \hspace{-0.02cm}\in\hspace{-0.02cm} X) \hspace{-0.02cm}= \hspace{-0.02cm} \begin{cases} 1 \hspace{-0.02cm}-\hspace{-0.02cm} \dfrac{k^*_C-\res}{(1-\overline{\varphi}_{x})S^{k^*}_C}, & \hspace{-0.15cm}\text{if } \hspace{-0.05cm} x \hspace{-0.02cm}\in\hspace{-0.02cm} \{x_1,\dots,x_{k^*_C}\}\\ 0, & \hspace{-0.15cm}\text{otherwise.}\end{cases} \label{sigma_marg}
\end{align}

 \end{theorem}

%
%

 Interestingly, from this theorem, we find that solving problem \eqref{eq:Q} indeed provides a lower bound on the optimal value of  \eqref{eq:P}. Furthermore, given a subset of sensor locations $C \subseteq \X$ of size greater than $r$, the post-security level that can be guaranteed for any component in $U_C$ is given by $1 - \frac{k^*_C-\res}{S^{k^*}_C}$. To achieve this post-security level, the network operator must ensure that the probability that each location in $C$ is monitored follows \eqref{sigma_marg}. Importantly, we find that the $k^*_C$ most critical locations must be monitored, with a probability that is decreasing with the associated lowest security level. In contrast, the remaining locations should not be monitored, as their security levels are already larger than the post-security level achieved for the $k^*_C$ most critical locations (by definition of $k^*_C$). We note that a monitoring strategy satisfying \eqref{sigma_marg} is guaranteed to exist by Farkas' lemma \cite{washburn1995two,2017arXiv170500349D}.
 
From Theorem \ref{prop:bounds}, we also find that upper bounds on the optimal value of the network monitoring problem \eqref{eq:P} can be computed from set packings. We will leverage this finding to efficiently compute optimality gaps for our solutions.

\subsection{Solution Approach}\label{sec:Approach}

From our structural analysis and Theorem \ref{prop:bounds}, we next derive a three-step approach to approximately solve the network monitoring problem \eqref{eq:P}. Steps 1 and 2 compute a monitoring strategy that optimally solves the relaxed problem \eqref{eq:Q}, and step 3 computes a set packing that achieves the upper bound from Theorem \ref{prop:bounds}.

In the first step, we formulate a MIP to compute $z^*_{(\Q)}$ as well as the marginal probabilities of monitoring each location at optimality of \eqref{eq:Q}. In particular, we leverage the fact that the objective function in \eqref{eq:R} can be expressed using these marginal probabilities. Indeed, for every $\sigma^1 \in \Delta^1$ and $x \in \X$,
\begin{align}
\mathbb{E}_{X\sim \sigma^1}[\bar{f}(X,x)] = \overline{\varphi}_x + (1- \overline{\varphi}_x)\mathbb{P}_{X\sim\sigma^1}(x \in X).\label{marg_f}
\end{align}


Using Lemma \ref{lem:Angle}, we then formulate the following MIP, which we refer to as a \emph{generalized covering set problem}:
%
%
%
%
%
%
%
%
%
%
\begin{alignat}{3}
\max_{y,z,\rho} \ & z \tag{$\mathcal{C}$}\label{eq:M}\\
\text{s.t.} \ &z \leq \overline{\varphi}_{x} + (1 - \overline{\varphi}_{x})\rho_x + M_1(1-y_x), \ &&\forall x \in \X\label{term1}\\
& z \leq \varphi_u + M_1\sum_{\quad \mathclap{\{x\in \X \, | \, u \in U_x\}}\quad } y_x, &&\forall u \in \U\label{term2}\\
& 0 \leq \rho_x \leq y_x, && \forall x \in \X\\
& \sum_{x \in \X} \rho_x = r,\label{resource_constraint}\\
& y_x \in \{0,1\}, && \forall x \in \X.\nonumber
\end{alignat}

Here, for every $x \in \X$, the binary variable $y_x$ determines if $x$ can receive a sensor, while the continuous variable $\rho_x$ determines the marginal probability that $x$ receives a sensor. $M_1$ is a ``big M'', which we can set to $1 - \min_{u \in \U} \varphi_u$. Given $C = \{x \in \X \ | \ y_x = 1\}$, we then observe that constraints  \eqref{term1} (combined with \eqref{marg_f}) model $\min_{x \in C}\mathbb{E}_{X\sim \sigma^1}[\bar{f}(X,x)]$, while constraints \eqref{term2} model $\min_{u \notin U_C} \varphi_u$. In fact, we show the equivalence between \eqref{eq:Q} and \eqref{eq:M}:




\begin{proposition}
\label{prop:LB} 
Problems \eqref{eq:Q} and \eqref{eq:M} have identical optimal values: $z^*_{(\Q)} = z^*_{(\mathcal{C})}$. Furthermore, a monitoring strategy $\sigma^{1^*} \in \Delta^1$ is an optimal solution of \eqref{eq:Q} if and only if there exists an optimal solution $(y^*,z^*,\rho^*)$ of \eqref{eq:M} satisfying $\rho^*_x = \mathbb{P}_{X\sim \sigma^{1^*}}(x \in X)$ for every $x \in \X$.
\end{proposition}

Thus, the optimal value of \eqref{eq:Q} can be computed by solving \eqref{eq:M}. In addition, given an optimal solution $(y^*,z^*,\rho^*)$ of \eqref{eq:M}, $\rho^*$ provides the marginal probabilities that each node must be monitored at optimality of \eqref{eq:Q}.

The second step of our approach consists of reconstructing an optimal solution of \eqref{eq:Q} from $\rho^*$. To this end, we utilize the combinatorial algorithm derived in \cite{dziubinski2018hide} that computes  a probability distribution $\sigma^{1^*} \in \Delta^1$ satisfying $\mathbb{P}_{X\sim \sigma^{1^*}}(x \in X) = \rho^*_x$ for every $x \in \X$. From Proposition \ref{prop:LB}, this guarantees the optimality of $\sigma^{1^*}$ for \eqref{eq:Q}.

Specifically, the algorithm iteratively expresses a vector $\rho \in [0,1]^{|\X|}$ satisfying \eqref{resource_constraint} as a convex combination of the incidence vector of a sensor positioning $X \in \A$, and another vector $\rho^\prime \in [0,1]^{|\X|}$ satisfying \eqref{resource_constraint} and with at least one additional integral component compared to $\rho^*$. We refer the reader to  \cite{dziubinski2018hide} for the detailed algorithm, which we henceforth refer to as the \emph{coordination algorithm}.


Importantly, given $n = |\supp(\rho^*)|$ the number of locations that must receive a sensor with positive probability, this algorithm runs in time $O(n^2)$ and constructs a monitoring strategy of support size $n+1$. We recall that a monitoring strategy with small support and/or small node basis can more easily be implemented in practice. We will numerically evaluate the node basis sizes of our solutions in our computational study.
%
%

The final step of our approach consists of evaluating the optimality gap of our solution for the network monitoring problem. To this end, we derive a MIP to efficiently compute the upper bound on $z^*_{(\mathcal{M})}$ from Theorem \ref{prop:bounds}. Specifically, we aim to compute a set packing $T \in \T$ that minimizes $1 - \frac{|T| - r}{S_T}$, which is equivalent to finding a set packing of size greater than $r$ that minimizes  $\frac{S_T}{|T| - \res}$. Although this is a \emph{nonlinear set packing problem}, we can reformulate it as the following MIP:
\begin{alignat}{3}
\min_{t,y,z}  & \displaystyle\sum_{l=\res+1}^{|\U|} \frac{1}{l - \res} t_l\tag{$\mathcal{P}$}\label{eq:N}\\
\text{s.t.}  & \sum_{u \in U_x} y_u \leq 1, && \forall x \in \X \label{cons_sp}\\
& \sum_{u \in \U} y_u = \sum_{l=\res + 1}^{|\U|} l z_l, \label{cons_nl_1}\\
& \sum_{l=\res + 1}^{|\U|} z_l = 1,\label{cons_nl_2}\\
&t_l \hspace{-0.05cm}\geq \hspace{-0.12cm}\sum_{u \in \U} (1\hspace{-0.05cm}-\hspace{-0.05cm}\varphi_u)^{-1} y_u \hspace{-0.05cm}-\hspace{-0.05cm} M_2 (1\hspace{-0.03cm}-\hspace{-0.03cm}z_l), && \forall l \hspace{-0.05cm}\in\hspace{-0.05cm} \llbracket \res \hspace{-0.05cm}+\hspace{-0.05cm} 1,|\U|\rrbracket\label{cons_product_1}\\
& t_l \geq 0, && \forall l \hspace{-0.05cm}\in\hspace{-0.05cm} \llbracket \res \hspace{-0.05cm}+\hspace{-0.05cm} 1,|\U|\rrbracket\label{cons_product_2}\\
& z_l \in \{0,1\}, &&  \forall l \hspace{-0.05cm}\in\hspace{-0.05cm} \llbracket \res \hspace{-0.05cm}+\hspace{-0.05cm} 1,|\U|\rrbracket\nonumber\\
& y_u \in \{0,1\}, && \forall u \in \U.\nonumber
\end{alignat}

Specifically, the binary variables $y_u$, for $u \in \U$ select the components to be part of the set packing. To parametrize the nonlinearity of the problem, we introduced for every $l \in  \llbracket \res + 1,|\U|\rrbracket$ the binary variable $z_l$, which is equal to 1 if the set packing is of size $l$. Finally, for every $l\in  \llbracket \res + 1,|\U|\rrbracket$, the continuous variable $t_l$ represents the product $z_l \sum_{u \in \U} (1-\varphi_u)^{-1}y_u$. Constraints \eqref{cons_sp} ensure that $\supp(y)$ is a set packing. Constraints \eqref{cons_nl_1}-\eqref{cons_nl_2} select the set packing size. Then, constraints \eqref{cons_product_1}-\eqref{cons_product_2} ensure that at optimality, the objective value is equal to ${\sum_{u \in \supp(y)} (1-\varphi_u)^{-1}}/{(|\supp(y)| - r)}$. $M_2$ is also a ``big M'', which we can set to $|\U| (1 - \max_{u \in \U}\varphi_u)^{-1}$. 

As a consequence, the upper bound on the optimal value of \eqref{eq:P} from Theorem \ref{prop:bounds} can computed by solving \eqref{eq:N}: 



%

\begin{proposition}
\label{prop:UB}
The following equality holds:
$$\min\bigg\{1,\min_{T \in \T} 1 -  \frac{|T|-\res}{S_T}\bigg\} = 1 - \frac{1}{z^*_{(\mathcal{P})}}.$$
\end{proposition}

We note that the number of variables in \eqref{eq:N} can be reduced in practice by computing an upper bound on the maximum set packing size in the networked system. Given the maximum set packing size $m^*$ (or an upper bound), then we can reduce the variables in \eqref{eq:N} from $|\U| + 2 (|\U|- r)$ to $|\U| + 2 (m^*-r)$, and set $M_2$ to $m^*(1 - \max_{u \in \U}\varphi_u)^{-1}$, which can provide significant computational benefits.

Thus, we obtain a three-strep approach to compute an approximate solution to the network monitoring problem \eqref{eq:P} and evaluate its optimality gap. Next, we further study our solution for problem instances with particular characteristics.


\subsection{Special Cases}
\label{subsec:exactNE}

\subsubsection{Non-overlapping monitoring sets}
We consider problem instances for which monitoring sets are pairwise disjoint, i.e., $U_x \cap U_{x^\prime} = \emptyset$ for every $x \neq x^\prime \in \X$. Such instances occur when it is desirable to reduce sensor interference or the energetic cost of the network \cite{cardei2005improving,rs6010740}. In other contexts such as in security games, disjoint monitoring is naturally satisfied \cite{powell2009sequential,behnezhad2018battlefields,musegaas2022stackelberg}.
For such instances, we show that our approach in Section~\ref{sec:Approach} optimally solves the network monitoring~problem:
\begin{proposition}\label{theorem:analytical_solution_special_case}
If the monitoring sets are pairwise disjoint, then $z^*_{(\Q)} = z^*_{(\M)}$. Furthermore, any monitoring strategy satisfying \eqref{sigma_marg} for $C = \X$ is an optimal solution of \eqref{eq:P}.
\end{proposition}

For instances when monitoring sets do not overlap, we find that problem \eqref{eq:P} can be optimally solved by first selecting the $k^*_\X$ most critical sensor locations and monitoring them with probabilities given by \eqref{sigma_marg}; and then by utilizing the coordination algorithm to construct a monitoring strategy that randomizes the placement of $r$ sensors and is consistent with the marginal probabilities in \eqref{sigma_marg}.

\subsubsection{Homogeneous security levels} We now consider instances where the security levels are identical across the components in $\U$. In such cases, we refine the bounds on the optimal value of \eqref{eq:P} using the minimum set cover and maximum set packing sizes ($n^*$ and $m^*$):
\begin{proposition}
\label{prop:identical}
If the security levels are identical and equal to $\varphi$, then:
\begin{align*}
z^*_{(\Q)} \hspace{-0.05cm}=\hspace{-0.05cm} \varphi + (1-\varphi) \frac{r}{n^*} \hspace{-0.05cm}\leq\hspace{-0.03cm} z^*_{(\M)} \hspace{-0.05cm} \leq \hspace{-0.03cm} \min\bigg\{1,\varphi + (1-\varphi) \frac{r}{m^*}\bigg\}.
\end{align*}
\end{proposition}

In such instances, we find that \eqref{eq:M} and \eqref{eq:N} can be simplified and the bounds on $z^*_{(\M)}$ can be computed by solving a minimum set cover and maximum set packing problem, as in \cite{2017arXiv170500349D}. If $C^*$ denotes a minimum set cover, then a monitoring strategy that optimizes $z^*_{(\Q)}$ monitors each location in $C^*$ with identical probability (given by \eqref{sigma_marg}), and can be constructed by cycling the $r$ sensors over $C^*$ (see \cite{2017arXiv170500349D} for reference).


%
%
%

  \section{Alternative Approaches}\label{sec:Alt}
  
  In order to compare our solution approach derived in Section~\ref{sec:Approach}, we next consider two alternative approaches for solving the network monitoring problem. In particular, we show how the column generation algorithm and multiplicative weights update algorithm can be applied to solve \eqref{eq:P}.
  

    \subsection{Column Generation Algorithm}
Problem \eqref{eq:P} can be formulated as a linear program with a large number (${|\X| \choose r} + 1$) of variables and a small number ($|\U| +1$) of constraints. Such problems can potentially be solved using the column generation algorithm, as long as the pricing problem can be solved efficiently.

Specifically, at each iteration of the column generation algorithm, the following restricted master problem is solved:
\begin{alignat}{3}
\max \ & z \tag{$\M_\mathcal{I}$}\label{eq:P_I}\\
\text{s.t.} \ & z - \sum_{X\in \mathcal{I}}\sigma^1_X f(X,u) \leq 0, \ && \forall u \in \U\label{dual1}\\
& \sum_{X \in \mathcal{I}}\sigma^1_X = 1,\label{dual2}\\
& \sigma^1_X \geq 0, && \forall X \in \mathcal{I},\nonumber
\end{alignat}
where $\mathcal{I} \subseteq \A$ is a subset of feasible sensor positionings. To determine if the optimal solution of \eqref{eq:P_I} is optimal for \eqref{eq:P}, the algorithm then determines the variable in \eqref{eq:P} with highest reduced cost. Let $(\alpha^*,\beta^*)$ denote the optimal dual variables of \eqref{eq:P_I}, associated with constraints \eqref{dual1}-\eqref{dual2}. Then, the variable with highest reduced cost can be determined by solving the following pricing problem:
\begin{align*}
\max_{X\in\A} - \beta^* +\sum_{u\in\U} \alpha^*_u f(X,u) .
\end{align*}
In fact, the pricing problem for \eqref{eq:P} is a maximum weighted covering problem, and can be solved using the following IP:
\begin{alignat}{3}
\max_{y,z} \ &  - \beta^* + \sum_{u \in \U} \alpha^*_u (\varphi_u  + (1-\varphi_u)y_u)\tag{$\mathcal{W}$}\label{eq:O}\\
\text{s.t.} \ & y_u \leq \sum_{\{x \in \X \, | \, u \in U_x\}} z_x,\ && \forall u \in \U\label{covering_pb_1}\\
& \sum_{x \in \X} z_x = r,\\
&y_u \in \{0,1\}, && \forall u \in \U\\
& z_x \in \{0,1\}, && \forall x \in \X\label{covering_pb_4}.
\end{alignat}

In this IP, the binary variables $z_x$, for $x \in \X$, represent the nodes that receive a sensor, and the binary variables $y_u$, for $u \in \U$, represent the components that are monitored by at least one sensor. 

If $z^*_{(\mathcal{W})} = 0$, then the optimal solution of \eqref{eq:P_I} is optimal for \eqref{eq:P}, and the algorithm terminates. If instead $z^*_{(\mathcal{W})} > 0$, then given an optimal solution $(y^*,z^*)$ of \eqref{eq:O}, the algorithm adds $X^* = \supp(z^*)$ to the set $\mathcal{I}$ of variable indices, and solves the new restricted master problem \eqref{eq:P_I}.

  \subsection{Multiplicative Weights Update Algorithm}

Problem \eqref{eq:P} is also equivalent to a simultaneous two-person zero-sum game, where the first player is the network operator who places $r$ sensors across the networked system and the second player is an attacker who targets a component. The payoff function is given by $f$, which the operator (resp. attacker) wants to maximize (resp. minimize).

In fact, \cite{FREUND199979} devised a multiplicative weights update algorithm for computing approximate mixed Nash equilibria of such simultaneous games, as long as one player has a small number of actions and the other player can efficiently compute best responses. 
%
The algorithm runs for $N$ iterations, and at each iteration $t \in \llbracket 1,N\rrbracket$, generates a probability distribution $\sigma^{2,t}$ over $\U$ for the attacker as well a feasible sensor positioning $X^t \in \A$ for the network operator. The algorithm starts by setting $\sigma^{2,1}$ to the uniform distribution over $\U$. Then, at each iteration $t \in \llbracket 1,N\rrbracket$, the algorithm determines the operator's best response to $\sigma^{2,t}$ by solving the following problem:
%
%
%
$$X^t \in \arg \max_{X \in \A} \mathbb{E}_{u \sim \sigma^{2,t}}[f(X,u)].$$
In our setting, this problem can also be formulated as a maximum weighted covering problem and can be solved using the following IP:
\begin{alignat}{3}
\max_{y,z} \ & \sum_{u \in \U} \sigma^{2,t}_u (\varphi_u  + (1-\varphi_u)y_u)\tag{$\mathcal{W}^\prime$}\label{eq:W'}\\
\text{s.t.} \ & \eqref{covering_pb_1}-\eqref{covering_pb_4}.\nonumber
\end{alignat}

Then, the algorithm updates the attacker's strategy using the following multiplicative weights update:
$$\forall u \in \U, \ \sigma^{2,t+1}_u = \sigma^{2,t}_u \frac{\eta^{f(X^t,u)}}{\sum_{u^\prime \in \U} \sigma^{2,t}_{u^\prime} \eta^{f(X^t,{u^\prime})}},$$
where $\eta = \left(1 + \sqrt{2 \frac{\ln |\U|}{N}}\right)^{-1}$. 
The authors in \cite{FREUND199979} showed that the monitoring strategy $\bar{\sigma}^1 \in \Delta^1$ that uniformly randomizes over the $N$ best responses $X^1,\dots,X^N$ satisfies:
\begin{align*}
\min_{u \in \U} \mathbb{E}_{X\sim \bar{\sigma}^1}[f(X,u)] \geq z^*_{(\mathcal{M})} - \sqrt{\frac{2 \ln |\U|}{N}} - \frac{\ln|\U|}{N}.
\end{align*}

Thus, if $N = 4\lceil \frac{\ln |\U|}{\epsilon^2}\rceil$ for $\epsilon > 0$, then the (absolute) optimality gap associated with the monitoring strategy $\bar{\sigma}^1$ generated by this algorithm is upper bounded by $\epsilon$.

In fact, \cite{krause2011randomized} showed that each iteration of this algorithm can be conducted in polynomial time by approximately solving \eqref{eq:W'} using a greedy algorithm. However, the resulting guarantee on the monitoring strategy $\bar{\sigma}^1$ becomes
\begin{align*}
\min_{u \in \U} \mathbb{E}_{X\sim \bar{\sigma}^1}[f(X,u)] \geq (1-1/e) z^*_{(\mathcal{M})} - \epsilon,
\end{align*}
for $T = 4\lceil \frac{\ln |\U|}{\epsilon^2}\rceil$. Since commercial solvers are now highly efficient at solving classical IPs, we instead optimally solve \eqref{eq:W'} at each iteration of the algorithm.

We now move to computational experiments, where we evaluate the three approaches presented in this article and compare their running times, as well as the performance and implementability of the monitoring strategies they generate.

\section{Computational Study}\label{section:computations}

In this section, we describe how our strategic network monitoring problem can be applied for contamination detection in water distribution systems and actuator protection in networked control systems. Our experiments are coded in Julia and all optimization problems are solved using the Gurobi solver v9.0 on a computer with 2.3-GHz, 8-Core Intel Core i9 processor and 32 GB of RAM.

\subsection{Application I: Contamination Detection in Water Distribution Systems}\label{section:simulations1}


We first consider the security problem where the operator of a water distribution system aims to allocate sensors to monitor their network against adversarially-induced contamination events. We model the water distribution system as a directed graph, with the vertices representing pumps, junctions, and water tanks, while the edges represent pipes. 
The direction of every edge is adopted to be in the direction of the water flow. 
We assume that the operator has access to contaminant detection sensors that can be deployed on valve access points or fire hydrants and measure water quality indicators such as electrical conductivity, free and total chlorine, turbidity, and oxygen reduction potential in the water  \cite{1219475,Fischer2017,AISOPOU2012235}.




In this application, the set of vulnerable components $\U$ (resp. set of sensor locations $\X$) is given by the set of edges (resp. nodes) of the graph modeling the water network. Monitoring sets are constructed through simulations using a hydraulic network solver \cite{epanet:2002} that tracks the advection and reaction dynamics from a contaminant intrusion event \cite{doi:10.1021/es3014024}.
%
%
%
%
In practice, the security levels of components can be assessed based on previously deployed security measures and their accessibility to an adversary. One can then use different security scales to determine quantitative values for the security levels based on that assessment \cite{RiskAssesment}. In this application, we adopt the scale $\{0.2,0.4,0.6,0.8\}$ for security levels, and randomly assign them to each component.



%

\begin{figure}[htbp]
 \begin{tikzpicture}[scale=1]
\begin{axis}[
  grid=major, 
  grid style={dashed,gray!30}, 
  xlabel=Number of sensors $r$,
  ylabel={Lowest post-security level $\min\limits_{u \in \U}\mathbb{E}_{X\sim \sigma^1}[f(X,u)]$},
   font=\small,
  xmin=0,
  xmax=19,
  ymin=0.2,
  ymax=1,
 ytick={0.2, 0.4,0.6,0.8,1},
  xtick={0, 2,4,6,8,10,12,14,16,18},
  legend pos=south east,
	legend cell align={left},
  legend style={draw=none, font = \footnotesize},
]

  \addplot[color = red, very thick] table[x=rs, y=Upper Bound, col sep=comma, comment chars={\%}] {ky5_results_obj.csv};

  
\addplot[color = black, very thick, dashed] table[x=rs, y=Optimal Value, col sep=comma, comment chars={\%}] {ky5_results_obj.csv};

  \addplot[color = blue, very thick] table[x=rs, y=Lower Bound Marg, col sep=comma, comment chars={\%}] {ky5_results_obj.csv};

 \legend{$z^*_{(\mathcal{P})}$,$z^*_{(\mathcal{M})}$, $z^*_{(\mathcal{C})}$}

 
\end{axis}

\end{tikzpicture}   
\caption{Illustration of the bounds from Theorem \ref{prop:bounds} for varying numbers of sensors.}
\label{fig:bounds}
\end{figure}
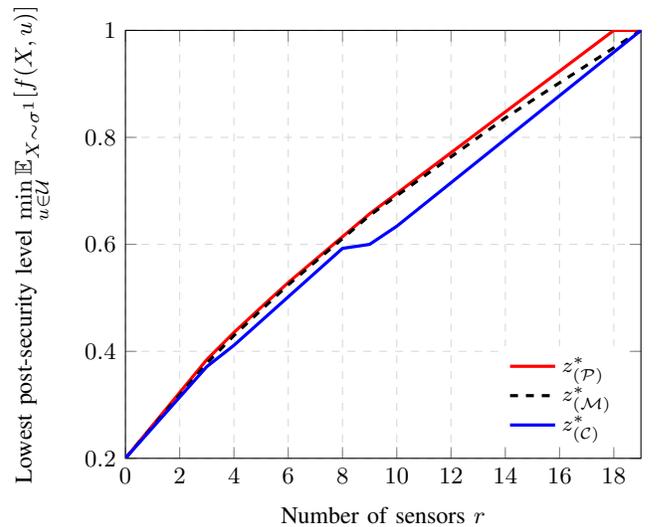

 Next, we implement the three solution approaches described in Sections \ref{sec:Method} and \ref{sec:Alt} on the water distribution system ky5~\cite{jolly2013research}. This anonymized real-world water network from Kentucky comprises 496 pipes and 420 nodes (junctions, water tanks, pumps), satisfies a demand of 2.28 million gallons of water per day, and spans 52.3 miles. In Fig. \ref{fig:bounds}, we illustrate the bounds derived in Theorem \ref{prop:bounds} and computed by solving the generalized covering set problem \eqref{eq:M} and the nonlinear set packing problem \eqref{eq:N}, respectively.

From Fig. \ref{fig:bounds}, we observe that the bounds obtained from \eqref{eq:M} and \eqref{eq:N} are close to the optimal value of \eqref{eq:P}. On average, the optimality gap given by the bounds does not exceed $4.7\,\%$, while the actual optimality gap is equal to $3.6\,\%$ on average.
 
 In general, we observe that $z^*_{(\M)}$ increases nonlinearly with respect to the number of sensors. When $r$ is small, the optimal monitoring strategy focuses on protecting the most vulnerable components (with a security level of $0.2$). As the operator has access to more sensors, they randomize their placement to increase the post-security level of less critical components. With 19 sensors, the operator can monitor all components and ensure a post-security level of 1 for the entire network.

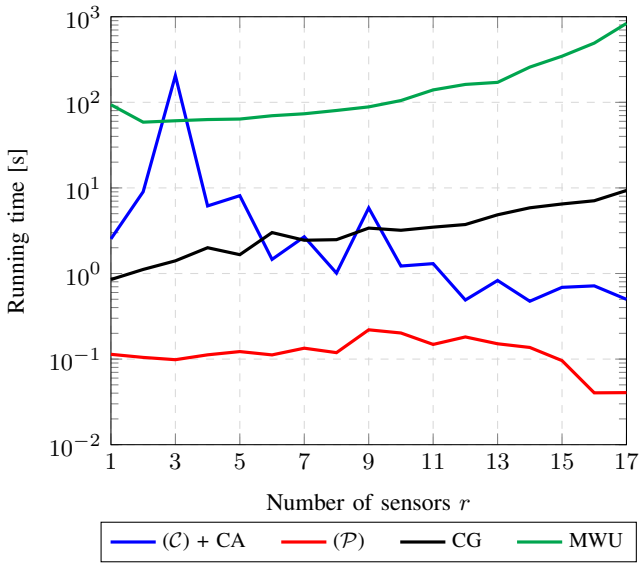
\begin{figure}[htbp]
 \begin{tikzpicture}[scale=1]
\begin{semilogyaxis}[
  grid=major, 
  grid style={dashed,gray!30}, 
  xlabel=Number of sensors $r$,
  ylabel={Running time [s]},
   font=\small,
  xmin=1,
  xmax=17,
  ymin=0.01,
  ymax=1000,
  xtick={1,3,5,7,9,11,13,15,17},
legend columns=4,
  legend style={/tikz/column 2/.style={
                column sep=0pt,
            }, font = \footnotesize,legend cell align={left},at={(0.5,-0.18)},anchor=north},
]

\legend{\eqref{eq:M} + CA \quad \quad, \eqref{eq:N} \quad \quad, CG \quad \quad,  MWU}

  \addplot[color = blue, very thick] table[x=rs, y=Time Lower Bound Marg, col sep=comma, comment chars={\%}] {ky5_results_times.csv};

  \addplot[color = red, very thick] table[x=rs, y=Time Upper Bound, col sep=comma, comment chars={\%}] {ky5_results_times.csv};

\addplot[color = black, very thick] table[x=rs, y=Time Optimal Value, col sep=comma, comment chars={\%}] {ky5_results_times.csv};

    \addplot[color = Green, very thick] table[x=rs, y=Time Multiplicative Weight Update, col sep=comma, comment chars={\%}] {ky5_results_times.csv};

\end{semilogyaxis}

\end{tikzpicture}   
\caption{Approaches' running times for varying numbers of sensors.}
\label{fig:times}
\end{figure}

Next, we compare the running times of the various approaches. Specifically, we illustrate the time to (a) solve \eqref{eq:N} and get an upper bound on $z^*_{(\M)}$, (b) solve \eqref{eq:M} and construct a monitoring strategy using the coordination algorithm (CA), (c) solve \eqref{eq:P} using column generation (CG), and (d) solve \eqref{eq:P} using multiplicative weights updates (MWU) with $N=2,480$ iterations to ensure a theoretical (absolute) optimality gap of $\epsilon = 0.1$. Their running times are illustrated in Fig. \ref{fig:times}.

%

Fig. \ref{fig:times} shows that the upper bound on $z^*_{(\M)}$ derived in Theorem \ref{prop:bounds} can be efficiently computed; by solving \eqref{eq:N}, we obtain the upper bound in less than $0.2$ seconds. We also observe that the number of sensors has limited impact on the running time. 

Furthermore, we find that CG is approximately 100 times faster than MWU. This is in part due to the number of iterations that are required for MWU to theoretically guarantee a small optimality gap. Interestingly, we find that CG can efficiently compute an optimal solution of \eqref{eq:P} when the number of sensors is small. However, as $r$ increases, the running time of CG increases exponentially. Finally, we find that our approximate monitoring strategy, which achieves the lower bound in Theorem \ref{prop:bounds} and is obtained by solving \eqref{eq:M} + CA, can be efficiently computed. 

\begin{figure}[htbp]
 \begin{tikzpicture}[scale=1]
\begin{axis}[
  grid=major, 
  grid style={dashed,gray!30}, 
  xlabel=Number of sensors $r$,
  ylabel={Node basis size $|\X(\sigma^1)|$},
   font=\small,
  xmin=1,
  xmax=17,
  ymin=0,
  ymax=200,
  xtick={1,3,5,7,9,11,13,15,17},
  legend style={font = \footnotesize,draw=none,legend cell align={left},at={(0.8,0.37)},anchor=north},
]

\legend{\eqref{eq:M} + CA, CG, MWU}


  \addplot[color = blue, very thick] table[x=rs, y=Node Basis Size Marg, col sep=comma, comment chars={\%}] {ky5_results.csv};
  
\addplot[color = black, very thick] table[x=rs, y=Node Basis Size at Optimality, col sep=comma, comment chars={\%}] {ky5_results.csv};

    \addplot[color = Green, very thick] table[x=rs, y=Node Basis Size Multiplicative Weight Update, col sep=comma, comment chars={\%}] {ky5_results.csv};

\end{axis}

\end{tikzpicture}
\caption{Node basis sizes for varying numbers of sensors.}
\label{fig:node_basis}
\end{figure}
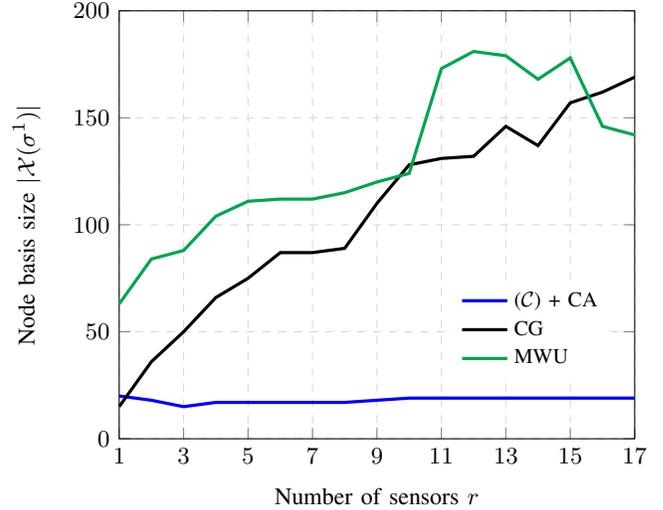

Finally, we compare the implementability of the solutions generated by the three approaches by illustrating their node basis sizes in Fig. \ref{fig:node_basis}.

Fig. \ref{fig:node_basis} clearly shows the implementability advantage of our monitoring strategy. In particular, we find that our solution requires monitoring at most 20 different locations and achieves an excellent worst-case post-security level (from Fig. \ref{fig:bounds}). In contrast, the optimal post-security level achieved using CG or MWU requires monitoring up to 180 different locations. From an implementation perspective, our monitoring strategy is considerably easier to translate into a daily/weekly scheduling of monitoring operations.


\subsection{Application II: Protecting Actuators from Extended Replay Attacks} \label{section:simulations2}

We now show how our findings can be used to protect actuators in a networked control system from extended replay attacks. The networked system consists of a set of physical states $\mathcal{X}$ and actuators $\mathcal{U}$, as well as fixed unprotected sensors.
We consider an attacker who aims to conduct an undetectable attack on an actuator using an extended replay strategy. More precisely an extended replay attack on an actuator is undetectable if the sensor measurements remain in the steady state during the attack. To this end, the attacker must also compromise a subset of components (i.e., unprotected sensors and other actuators) to ensure that the attack is undetected. Such subset of components typically depends on the dynamics of the networked system, and can be efficiently determined in large-scale systems \cite[Chapter 7]{milovsevic2020security},\cite{TAC_2019_Jezdimir}.

Importantly, the vulnerability of an actuator $u$ with respect to extended replay attacks can be determined based on its security index $\delta_u$, which represents the minimum number of sensors and actuators, including $u$ itself, that the attacker must compromise to conduct an undetectable extended replay attack against $u$, in the absence of protected sensors \cite[Chapter 7]{milovsevic2020security}. If no undetectable extended replay attack against $u$ is possible, then $\delta_u = +\infty$. We can then scale these security indices and define security levels for this application, following our modeling framework. Specifically, we adopt an identical scale as in the previous application, and define the security level of each actuator $u \in \U$ as follows:
 \begin{align*}
\varphi_u =&\textstyle 0.2 \cdot \mathds{1}_{[\delta_u \in (0,5] ]}  
+ 0.4 \cdot \mathds{1}_{\{\delta_u \in (5,15] \}}+ 0.6 \cdot \mathds{1}_{\{\delta_u \in (15,20] \}}   \\
&\textstyle
+ 0.8 \cdot \mathds{1}_{\{\delta_u \in (20,+\infty)\}}+ 1.0 \cdot \mathds{1}_{\{\delta_u= +\infty\}}.
\end{align*}

\begin{table*}[!htbp]
\caption{Comparison of solution approaches on IEEE 2383 bus power system.}
\label{tab:ieee}
\centering
\begin{tabular}{ccccccccccccc}
\hline
{\# sensors } &  & \multicolumn{3}{c}{Optimality gap [\%]} &  & \multicolumn{3}{c}{Running time [s]} &  & \multicolumn{3}{c}{Node basis size} \\ \cline{3-5} \cline{7-9} \cline{11-13} 
 $r$ &  & \eqref{eq:M} + CA & CG  & MWU &  & \eqref{eq:M} + CA & CG & MWU &  & \eqref{eq:M} + CA & CG & MWU \\ \hline
\phantom{1}50 &  &0.0 & \phantom{1}0.0 & 0.3  &  & \phantom{11}1.6 & 71.3 & 227.7 &  & 223 & 373 & \phantom{1}430 \\
100 &  &0.0 & \phantom{1}0.4 & 0.4  &  & \phantom{11}1.9 & 600 & 245.1 &  & 318 & 868 & \phantom{1}737 \\
150 &  &0.0 & \phantom{1}1.1 & 0.4  &  & \phantom{1}13.2 & 600 & 302.6 &  & 328 & 981 & \phantom{1}824 \\
200 &  &0.0 & \phantom{1}8.0 & 0.3  &  & \phantom{11}5.1 & 600 & 412.2 &  & 328 & 1009 & \phantom{1}827 \\
250 &  &0.0 & 15.0 & 0.3  &  & 174.2 & 600 & 581.2 &  & 344 & 1024 & 1115 \\
300 &  &0.0 & 24.0 & 0.3  &  & \phantom{1}39.2 & 600 & 591.0 &  & 344 & 1086 & \phantom{1}612 \\ \hline
\end{tabular}
\end{table*}

To defend the networked system, we assume the operator possesses protected sensors that cannot be compromised by the attacker and that can be positioned within the network. We assume that each positioned sensor can only measure one physical state, which is a commonly adopted model of large--scale systems~\cite{7122316,6762966,7524914}. In particular, a protected sensor can detect if the state it measures deviates from its steady state value during an extended replay attack. Thus, actuator $u \in \U$ belongs to the monitoring set of state $x \in \X$ if placing a protected sensor at $x$ prevents the attacker from conducting an undetectable extended replay attack on $u$.

We now consider the IEEE 2383 bus power system.  
We model the system using linearized swing equations where the generators are represented by two states (rotor angle and frequency) and load buses with one state (voltage angle)~\cite{4110445}. 
We assume that all $|\X| = 3037$ states are measurable, and that the attacker can conduct an attack using $|\U|=1042$ loads (we randomly selected 30$\%$ of the loads to be attackable) \cite{5976424}.  
Table \ref{tab:ieee} compares the performance of the different solution approaches on this networked system. Note that we set a time limit of 10 minutes for each approach.

We observe that for this network, our solution approach optimally solves the network monitoring problem \eqref{eq:P}. Indeed, the upper bound obtained by solving \eqref{eq:N} guarantees that the monitoring strategy computed from \eqref{eq:M} + CA is optimal for \eqref{eq:P} for every number of sensors. Although we observe some variance in the running time of our approach, the average running time for this network is 29 seconds. In general, we find that the hardest instances to solve are for which the number of sensors to position is high but not too close to the minimum set cover size.
%
%
%
%
On the other hand, MWU achieves a very high performance, but requires 393 seconds on average. This illustrates the fact that the empirical performance of MWU is better than the theoretical guarantees \cite{FREUND199979}. In contrast, we find that CG scales poorly with the number of sensors, which is interesting given that MWU and CG solve the very similar IPs \eqref{eq:O} and \eqref{eq:W'}. Finally, we similarly note that the node basis (i.e., the locations that are monitored with positive probability) of our monitoring strategy is up to 3 times smaller than when utilizing CG or MWU, thus providing a solution that is more easily implementable in practice.

%
%

\subsection{Summary and Extension of Solution Methodology}

In summary, we find that by extracting the structure of the network monitoring problem, we could derive an efficient solution approach. Our approach solves two IPs and runs the coordination algorithm to directly obtain a monitoring strategy with strong performance guarantees and that is easier to implement in practice. This contrasts with the classical CG and MWU algorithms that iteratively refine their solutions and generally run more slowly.

Interestingly, more flexibility can be embedded in our solution approach. For instance, if the operator is interested in an even simpler monitoring strategy, we can add the cardinality constraint $\sum_{x \in \X} y_x \leq s$ in \eqref{eq:M} to find the best monitoring strategy of \eqref{eq:Q} that monitors at most $s$ locations. This permits the operator to achieve the desired tradeoff between solution implementability and guaranteed post-security level.

Finally, we note that the three solution approaches can be combined for an even better performance. Indeed, for the instances for which our approach does not provide an optimal solution to the network monitoring problem \eqref{eq:P}, we can utilize our solution as a warm start for CG and/or MWU, which will then iteratively improve it.

\section{Conclusion}  \label{section:ch8_conclusion}

In this article, we studied the problem of strategically monitoring large-scale networked systems with heterogeneous component criticalities using multiple sensors. Specifically, we formulated a large-scale maximin problem where the network operator selects a randomized placement of their protected sensors to improve the post-security level of the network's most critical components against an attacker. 

We developed a three-step solution approach that leverages the structure of the network monitoring problem. First, we solved a generalized covering set problem to determine the marginal probabilities to monitor each location. Then, we utilized a combinatorial algorithm to randomize the positioning of the multiple sensors and obtain a monitoring strategy that is consistent with the marginal probabilities. Finally, we solved a nonlinear set packing problem to evaluate our solution's optimality gap. We also investigated problem instances with special structures, for which we refined our solution approach and optimality gaps.

 
 We then adapted two classical solution methods for solving our monitoring problem based on column generation and multiplicative weights updates, and compared them on real-world benchmark water distribution and power systems. Our computational study revealed that our solution approach outperforms the classical methods as the size of the problem increases. In particular, we find that our approach can efficiently compute near optimal solutions for large instances and are significantly simpler to implement. 


As part of future work, it would be interesting to investigate how the solution method developed in this article (combining smaller IPs and combinatorial algorithms) can be applied to solve problems where an attacker targets multiple components simultaneously, or where sensors are faulty.


\appendix


\proof[Proof of Lemma \ref{lem:Angle}] We consider an optimal solution $\sigma^{1^*}$ of \eqref{eq:Q} and we set $C = \X(\sigma^{1^*})$. Then, $z^*_{(\mathcal{R}_C)} \geq \min_{x \in \X(\sigma^{1^*})}\mathbb{E}_{X\sim \sigma^{1^*}}[\bar{f}(X,x)]$. Furthermore, $\min_{u \notin U_C} \varphi_u = \min_{u \notin U_{\X(\sigma^{1^*})}} \varphi_u$. Therefore, $z^*_{(\Q)} \leq  \max_{C \subseteq \X} \min\{z^*_{(\mathcal{R}_C)},\min_{u \notin U_C} \varphi_u\}$.

To show the reverse inequality, let $C^* \in \arg\max_{C \subseteq \X} \min\{z^*_{(\mathcal{R}_C)},\min_{u \notin U_C} \varphi_u\}$. We also consider an optimal solution $\sigma^{1^*}$ of $(\mathcal{R}_{C^*})$. Without loss of generality, we can assume that $\X(\sigma^{1^*}) \subseteq C^*$ (i.e., $\sigma^{1^*}$ only places sensors within $C^*$).
%
%
%
 Then, $\min_{x \in \X(\sigma^{1^*})} \mathbb{E}_{X\sim\sigma^{1^*}}[\bar{f}(X,x)] \geq z^*_{(\mathcal{R}_{C^*})}$. Furthermore, let $u^\prime \notin U_{\X(\sigma^{1^*})}$. If $u^\prime \in U_{C^*}\setminus U_{\X(\sigma^{1^*})}$, let $x^\prime \in C^*\setminus\X(\sigma^{1^*})$ such that $u^\prime \in U_{x^\prime}$. We then obtain:
\begin{align*}
\varphi_{u^\prime} \geq \overline{\varphi}_{x^\prime} = \mathbb{E}_{X \sim \sigma^{1^*}}[\bar{f}(X,x^\prime)] \geq z^*_{(\mathcal{R}_{C^*})}.
\end{align*}

If instead $u^\prime \notin U_{C^*}$, then $\varphi_{u^\prime} \geq \min_{u \notin U_{C^*}} \varphi_u$. 
In conclusion, $z^*_{(\Q)} =  \max_{C \subseteq \X} \min\{z^*_{(\mathcal{R}_C)},\min_{u \notin U_C} \varphi_u\}$.
\endproof

 \proof[Proof of Theorem \ref{prop:bounds}] We first show that $z_{(\M)}^* \geq z_{(\Q)}^*$: Let $\sigma^{1^*} \in \Delta^1$ be an optimal solution of \eqref{eq:Q} and let $u^\prime \in \U$. If $u^\prime \in U_{\X(\sigma^{1^*})}$, then let $x^\prime \in \X(\sigma^{1^*})$ such that $u^\prime \in U_{x^\prime}$. Since $\varphi_{u^\prime} \geq \overline{\varphi}_{x^\prime}$ and for every $X \in \A$, $ \mathds{1}_{\{x^\prime \notin X\}} \geq \mathds{1}_{\{u^\prime \notin U_X\}}$, then
\begin{align*}
 \mathbb{E}_{X\sim\sigma^{1^*}}[f(X,u^\prime)] &\geq \mathbb{E}_{X\sim\sigma^{1^*}}[1 - (1 - \overline{\varphi}_{x^\prime})  \mathds{1}_{\{x^\prime \notin X\}}]\\
& \geq \min_{\quad\mathclap{x \in \X(\sigma^{1^*})}\quad} \mathbb{E}_{X\sim\sigma^{1^*}}[1 - (1 - \overline{\varphi}_{x})  \mathds{1}_{\{x \notin X\}}].
\end{align*}

If instead $u^\prime  \notin U_{\X(\sigma^{1^*})}$, then:
\begin{align*}
&\mathbb{E}_{X\sim\sigma^{1^*}}[f(X,u^\prime)] = \varphi_{u^\prime} \geq \min_{u \notin U_{\X(\sigma^{1^*})} } \varphi_u.
\end{align*}

Thus, $z^*_{(\mathcal{M})} \geq z^*_{(\Q)}$. We next show that for every $C \subseteq \X$ of size greater than $r$, $z^*_{(\mathcal{R}_C)} = 1 - \frac{k^*_C-\res}{S^{k^*}_C}$.

We rewrite $C = \{x_1,\dots,x_n\}$ with $n > r$ and indices satisfying $\overline{\varphi}_{x_1} \leq \cdots \leq \overline{\varphi}_{x_n}$. Then necessarily, $k_C^* \geq  r +1$ since
\begin{align*}
1 - \frac{r + 1 - r}{\sum_{l=1}^{r+1} (1 - \overline{\varphi}_{x_l})^{-1}} \geq 1 - (1 - \overline{\varphi}_{x_{r+1}}) = \overline{\varphi}_{x_{r+1}}.
\end{align*}
 
Let $\sigma^{1^\prime} \in \Delta^1$ be a monitoring strategy satisfying \eqref{sigma_marg}.
%
%
 The existence of such a monitoring strategy follows from Farkas' lemma (for instance, see \cite[Lemma EC.4]{2017arXiv170500349D}), since
 \begin{align*}
 \sum_{l=1}^{k^*_C} 1 - \frac{k^*_C-\res}{(1-\overline{\varphi}_{x})S^{k^*}_C} = r,
 \end{align*}
 and for every $l \in \llbracket 1,k^*_C\rrbracket$,
 \begin{align*}
1 \geq 1 - \frac{k^*_C-\res}{(1-\overline{\varphi}_{x_l})S^{k^*}_C} \geq 1 - \frac{k^*_C-\res}{(1-\overline{\varphi}_{x_{k^*_C}})S^{k^*}_C} \geq 0.
 \end{align*}

We next consider $x \in C$. If $x \in \{x_1,\dots,x_{k^*_C}\}$, then:
\begin{align*}
\mathbb{E}_{X \sim \sigma^{1^\prime}}[\bar{f}(X,x)] &= \overline{\varphi}_x + (1-\overline{\varphi}_x)\mathbb{P}_{X \sim \sigma^{1^\prime}}(x \in X) \\
&= 1 - \frac{k^*_C-\res}{S^{k^*}_C}.
\end{align*}

If instead $x \in C\setminus\{x_1,\dots,x_{k^*_C}\}$, then:
\begin{align*}
\mathbb{E}_{X \sim \sigma^{1^\prime}}[\bar{f}(X,x)] &= \overline{\varphi}_x \geq \overline{\varphi}_{x_{k^*_C+1}} \\
&> 1 - \frac{k^*_C + 1 - r}{S^{k^*}_C + (1 - \overline{\varphi}_{x_{k^*_C+1}})^{-1}},
\end{align*}
which is equivalent to
\begin{align*}
\mathbb{E}_{X \sim \sigma^{1^\prime}}[\bar{f}(X,x)] > 1 - \frac{k^*_C - r}{S^{k^*}_C}.
\end{align*}

Thus, $z^*_{(\mathcal{R}_C)} \geq 1 - \frac{k^*_C-\res}{S^{k^*}_C}$.


Now, let us assume by contradiction that $z^*_{(\mathcal{R}_C)}> 1 - \frac{k^*_C-\res}{S^{k^*}_C}$. Let $\sigma^{1^*}$ be an optimal solution of $z^*_{(\mathcal{R}_C)}$. Then, we deduce that:
\begin{align*}
\forall x \in C, \ \mathbb{P}_{X\sim \sigma^{1^*}}(x \in X) > 1 - \frac{k^*_C-\res}{(1-\overline{\varphi}_x)S^{k^*}_C},
\end{align*}
 which induces the following contradiction:
 \begin{align*}
r &= \sum_{X \in \A} \sigma^{1^*}_X |X| = \sum_{x \in \X} \sum_{X \in \A} \sigma^{1^*}_X \mathds{1}_{\{x \in X\}} \\
&=\sum_{x \in \X} \mathbb{P}_{X\sim \sigma^{1^*}}(x \in X) > \sum_{l=1}^{k^*_C} 1 - \frac{k^*_C-\res}{(1-\overline{\varphi}_{x_l})S^{k^*}_C} = r.
\end{align*}

In conclusion, $z^*_{(\mathcal{R}_C)}= 1 - \frac{k^*_C-\res}{S^{k^*}_C}$, and an optimal solution is given by $\sigma^{1^\prime}$ satisfying \eqref{sigma_marg}.

We note that for every $C \subseteq \X$ of size no more than $r$, $z^*_{(\mathcal{R}_C)} = 1$ (which is obtained by placing a sensor on each location in $C$). Also note that in such cases, $\U\setminus U_C \neq \emptyset$ and $1 - \frac{k^*_C-\res}{S^{k^*}_C}\geq 1 \geq \min_{u \notin U_C} \varphi_u$. Thus, from Lemma \ref{lem:Angle}, we deduce that 
 $$z^*_{(\overline{\mathcal{M}})} = \max_{C \subseteq \X} \min\bigg\{1 - \frac{k^*_C-\res}{S^{k^*}_C},\min_{\mathclap{u \notin U_C}}\varphi_u\bigg\}.$$
 
 Finally, we prove the upper bound on $z^*_{(\M)}$. Since $f(X,u) \leq 1$ for every $X \in \A$ and every $u \in \U$, then $z^*_{(\M)} \leq 1$. Furthermore, by strong duality, we know that $z^*_{(\M)} = \min_{\sigma^2 \in \Delta^2} \max_{X \in \A}\mathbb{E}_{u\sim \sigma^2}[f(X,u)]$, where $\Delta^2 \coloneqq \{\sigma^2 \in [0,1]^{|\U|} \ | \ \sum_{u \in \U}\sigma^2_u = 1\}.$ Then, given a set packing $T \in \T$, we define $\sigma^{2^\prime} \in \Delta^2$ as follows:
 \begin{align*}
\sigma^{2^\prime}_u = \begin{cases} \frac{1}{(1-\varphi_u)S_T}, & \text{ if } u \in T\\
0, & \text{ otherwise.}\end{cases}
 \end{align*}
 
 We can then derive the following upper bound.
  \begin{align*}
  z^*_{(\M)} &\leq \max_{X\in \A}\mathbb{E}_{u\sim \sigma^{2^\prime}}[f(X,u)] \\
  &= \max_{X\in \A}\sum_{u \in T} \sigma^{2^\prime}_u (1 - (1-\varphi_u)\mathds{1}_{\{u \notin U_X\}})\\
  & = 1 - \frac{1}{S_T}\min_{X\in \A}\sum_{u \in T}\mathds{1}_{\{u \notin U_X\}} \leq 1 - \frac{|T|-r}{S_T}
 \end{align*}
 
In conclusion, $z_{(\M)}^*\leq \min\{1,\min_{T \in \T} 1 -  \frac{|T|-\res}{S_T}\}$.
 \endproof

\proof[Proof of Proposition \ref{prop:LB}] Let $(y^*,z^*,\rho^*)$ be an optimal solution of \eqref{eq:M}, and we set $C = \supp(y^*)$. Then, \eqref{eq:M} can be rewritten as follows:
\begin{alignat}{3}
\max_y &\min\bigg\{\max_\rho \min_{x \in \supp(y)} && \overline{\varphi}_x +  (1- \overline{\varphi}_x) \rho_x, \min_{u \notin U_{\supp(y)}} \varphi_u\bigg\}.\nonumber\\
 \text{s.t.} \ & 0 \leq \rho_x \leq y_x, \ &&\forall x \in \X\label{cons_1}\\
 & \sum_{x \in \X}\rho_x = r,\label{cons_2}\\
 & y_x \in \{0,1\}, \ &&\forall x \in \X.\nonumber
\end{alignat}
The equivalence between \eqref{eq:Q} and \eqref{eq:M} then follows from Lemma \ref{lem:Angle}, equality \eqref{marg_f}, and the existence of monitoring strategies satisfying constraints \eqref{cons_1}-\eqref{cons_2}.
\endproof

  \proof[Proof of Proposition \ref{prop:UB}]
We argued in the main text that Problem \eqref{eq:N} is indeed a MIP formulation of the following problem $\min_{T \in \T} \frac{S_T}{|T|-r}$. In particular, at optimality, constraints \eqref{cons_nl_2} and $M_2$ ensure that at optimality all $t_l$ variables are equal to 0 except for the index corresponding to the set packing size, thus ensuring the validity of the formulation. Then, $\min\{1,\min_{T \in \T} 1 -  \frac{|T|-\res}{S_T}\} = 1 - \frac{1}{z^*_{(\mathcal{P})}}.$ Note that we used the convention that if a \eqref{eq:N} is infeasible (which occurs if there is no set packing of size greater than $r$), then $z^*_{(\mathcal{P})} = + \infty$.
\endproof

\proof[Proof of Proposition \ref{theorem:analytical_solution_special_case}] We consider an instance where $U_x \cap U_{x^\prime} = \emptyset$ for every $x \neq x^\prime \in \X$. We index the locations in $\X = \{x_1,\dots,x_n\}$ in nondecreasing order of their criticalities. Then, we consider the set packing $T^\prime = \{u_1,\dots,u_{k^*_\X}\}\in \T$, where for every $l \in \llbracket 1,k^*_\X\rrbracket$, $\varphi_{u_l} = \overline{\varphi}_{x_l}$. In other words, $T^\prime$ selects the most critical component within each of the first $k^*_\X$ monitoring sets. As a consequence, $|T^\prime| = k^*_\X$ and $S_{T^\prime} = S^{k^*}_\X$. Since $U_{\X} = \U$, then from Theorem \ref{prop:bounds}, we obtain:
\begin{align*}
z^*_{(\mathcal{R}_\X)}  &= 1 - \frac{k^*_\X - r}{S^{k^*}_\X} \leq  z_{(\Q)}^*\leq z^*_{(\M)} \\
&\leq 1 - \frac{|T^\prime| - r}{S_{T^\prime}} = 1 - \frac{k^*_\X - r}{S^{k^*}_\X}.
\end{align*}

Thus, $z^*_{(\M)} = z^*_{(\Q)}$. Furthermore, any monitoring strategy $\sigma^{1^\prime} \in \Delta^1$ satisfying \eqref{sigma_marg} for $C = \X$, which is optimal for $(\mathcal{R}_\X)$, is then optimal for \eqref{eq:P}.
\endproof

\proof[Proof of Proposition \ref{prop:identical}]
We suppose that for every component, the security levels are identical and equal to $\varphi$. Let $C\subseteq \X$. Then, we find that $k^*_C = |C|$ since for every $k \in \llbracket 1,|C| \rrbracket$, $\frac{k-r}{k}(1-\varphi) \leq 1- \varphi$. Therefore, $z^*_{(\Q)} = \max_{C \subseteq \X} \min\{1 - \frac{|C| - r}{|C|}(1-\varphi), \min_{u \notin U_C} \varphi \}$.

Note that for any $C\subseteq \X$, $1 - \frac{|C| - r}{|C|}(1-\varphi) = \varphi + \frac{r}{|C|} (1-\varphi) \geq \varphi$. Therefore, 
\begin{align*}
z^*_{(\Q)} = \max_{C\in \mathcal{S}}  \varphi + \frac{r}{|C|} (1-\varphi) =  \varphi + \frac{r}{n^*}(1-\varphi),
\end{align*}
where $\mathcal{S}$ is the collection of set covers. Similarly, 
\begin{align*}
\min_{T\in\T} 1 - \frac{|T| - r}{S_T}  = \min_{T\in\T}  1 + \frac{r}{|T|}(1-\varphi) =  1 + \frac{r}{m^*}(1-\varphi).
\end{align*}

\endproof

\bibliographystyle{IEEEtran}

\bibliography{autosam} 

\vspace{-1.3cm}
\begin{biographynophoto}{Jezdimir Milo\v{s}evi\'{c}} received the M.Sc. degree in electrical engineering and computer science from the University of Belgrade, Belgrade, Serbia, in 2015 and the Ph.D. in electrical engineering and computer science from the KTH Royal Institute of Technology, Stockholm, Sweden, in 2020.

He is currently a research and development engineer in Scania's Autonomous Transport Solutions, Södertälje, Sweden. 
His research interests are within cyber-security of industrial control systems.
  \end{biographynophoto}

\vspace{-1.3cm}

\begin{biographynophoto}{Mathieu Dahan} received the M.S. and Ph.D. degrees in computational science and engineering from the Massachusetts Institute of Technology, Cambridge, MA, USA, in 2016 and 2019, respectively.

He is currently an Assistant Professor in the School of Industrial and Systems Engineering at the Georgia Institute of Technology, Atlanta, GA, USA. His research interests are in combinatorial optimization, game theory, and predictive analytics, with applications to service and healthcare operations, humanitarian systems, and logistics and supply chain management.

 \end{biographynophoto}

\vspace{-1.3cm}
 
\begin{biographynophoto}{Saurabh Amin} received the Ph.D. degree in systems engineering from the University of California, Berkeley, Berkeley, CA, USA, in 2011.

He is currently a Professor with the Department of Civil and Environmental Engineering. He is a member of the Laboratory for Information and Decision Systems with the Massachusetts Institute of Technology, Cambridge, MA, USA. His fields of expertise include control and optimization, applied game theory, and networks. His research focuses on the design and implementation of resilient monitoring and control algorithms for networked infrastructure systems.
 \end{biographynophoto}

\vspace{-1.3cm}

\begin{biographynophoto}{Henrik Sandberg} (F’23) received the M.Sc. degree in engineering physics and the Ph.D. degree in automatic control from Lund University, Lund, Sweden, in 1999 and 2004, respectively. 

He is currently a Professor with the Division of Decision and Control Systems, KTH Royal Institute of Technology, Stockholm, Sweden. His current research interests include security of cyber-physical systems, power systems, model reduction, and fundamental limitations in
control. 

Dr. Sandberg was a recipient of the Best Student Paper Award
from the IEEE Conference on Decision and Control in 2004, an Ingvar
Carlsson Award from the Swedish Foundation for Strategic Research in
2007, and a Consolidator Grant from the Swedish Research Council in
2016. He has served on the editorial boards of IEEE Transactions on Automatic Control and the IFAC Journal Automatica. and the IFAC Journal Automatica.
 \end{biographynophoto}

\end{document}